\documentclass{article}

\usepackage[a4paper]{geometry}

\usepackage{amsmath,amssymb}

\usepackage{setspace}

\usepackage{graphicx}

\newtheorem{lemma}{Lemma}
\newtheorem{theorem}{Theorem}
\newtheorem{definition}{Definition}

\numberwithin{equation}{section}

\title{\textbf{Discrete LIBOR Market Model Analogy}\footnotemark\footnote{This research has been done at Vienna University of Technology and Dublin City University. The author gratefully acknowledges the Austrian Christian Doppler Society (CD-laboratory PRisMa), as well as Science Foundation Ireland (Edgeworth Center and FMC2) for their support. Special thanks go to Dr. Friedrich Hubalek for critical discussion and advice on the subject.}}

\author{\small{\textbf{Andreas Hula}}\\ \footnotesize{School Of Mathematics, Dublin City University,Dublin 9,Dublin,Ireland}\\ \footnotesize{(andreas.hula@dcu.ie)}}

\date{\today}

\begin{document}

\doublespacing

\pagestyle{headings}

\maketitle

\tableofcontents
\vspace{0.6 cm}
\textbf{Mathematics Subjects Classification: 91G30, 60J75, 60G57}
\vspace{0.2 cm}
\begin{abstract}
 This paper provides a discrete time LIBOR analogue, which can be used for arbitrage-free discretization of L\'evy LIBOR models or discrete approximation of continuous time LIBOR market models. Using the work of Eberlein and {\"O}zkan \cite{EO2005} as an inspiration, we build a discrete forward LIBOR market model by starting with a discrete exponential martingale. We take this pure jump process and calculate the appropriate measure change between the forward measures.\\
Next we prove weak convergence of the discrete analogue to the continuous time LIBOR model, provided the driving process converges weakly to the continuous time one and the driving processes are PII's. \\
We also demonstrate the results of a practical implementation of this model and compare them to an implementation of an arbitrage free discretization by Glasserman and Zhao $\cite{GZ2000}$.
\end{abstract}
\section{Introduction}
LIBOR Market Models are intensely used in banking and finance since their appearance in the paper by Brace, Gatarek and Musiela $\cite{BGM}$. A rich literature has developed in the field, especially the book by Brigo and Mercurio $\cite{BM2006}$ should be mentioned as a accessible reference. Since all practical calculations need discretization of some sort, arbitrage free discretization as in Glasserman and Zhao $\cite{GZ2000}$ is desirable. With the papers by Jamshidian $\cite{J1999}$ and Eberlein and {\"O}zkan $\cite{EO2005}$ more general driving processes became the center of attention. The dissertation of Kluge $\cite{KLUGE}$ offers an insight into application and calibration of those models. However in continuous time, one is usually forced to carry out the so called "frozen drift approximation". In a time discrete model driven by discrete random variables, such an approximation could possibly be avoided, which was the original motivation to construct such a model.\\
It thus seemed natural to look for a fairly general discrete time analogue of these models. It should be arbitrage free and the analogy to the continuous time case should be well apparent. We treat this in section $(\ref{properties})$.\\
 Then we devote an equally important section to the convergence of a sequence of discrete analogues, under certain integrability and convergence assumptions, to a continuous time model. This property has been the starting point of our interest in discrete LIBOR models. We get a convergence result and can demonstrate, that our discrete model can be applied to find approximations for the exact (non frozen drift) L\'evy LIBOR model.\\
 The next section $(\ref{merits})$ gives an overview on the results from a model implementation and compares the results to a model based on the work $\cite{GZ2000}$. This is carried out in more detail in section $(\ref{implementation})$.\\
 The final section is devoted to an outlook on further possibilities for work based on these results.\\
 All proofs can be found in the appendix $(\ref{proofs})$.\\
 Apart from the papers above, the methods we use can be found mostly in the books by Protter $\cite{PROTTER}$ and Jacod and Shirjaev $\cite{JS2002}$. The theorems relevant for the convergence result can be found in appendix $(\ref{weak})$.\\
 If the readers interest is purely in the implementation results, I recommend to look at sections $3.2$, $4.1$ and $5$. If the readers interest is more in the convergence result, sections $3$ and $4$ cover this entirely.
 \section{The Merits Of The Approach}\label{merits}
 We can show that there is a discrete LIBOR market model analogue such that all LIBOR rates converge under their respective forward measures and also jointly under the terminal measure to a continuous time model, as soon as the driving process converges.\\
 The implementation of such a model gains us the possibility to calculate without the frozen drift approximation, assuming the driving process is simple enough (see section $(\ref{implementation})$ for details).\\
 Below you can see the results of $2$ implementations of this approach, compared to an implementation of one arbitrage free discretization by Glasserman and Zhao from $\cite{GZ2000}$ and implied volatility smiles derived from the data. The $3$ models have been chosen, such that all $3$ will be arbitrage free and converge to the lognormal LIBOR market model, without frozen drift.\\
 \includegraphics[width=4.8cm,height=4.5cm]{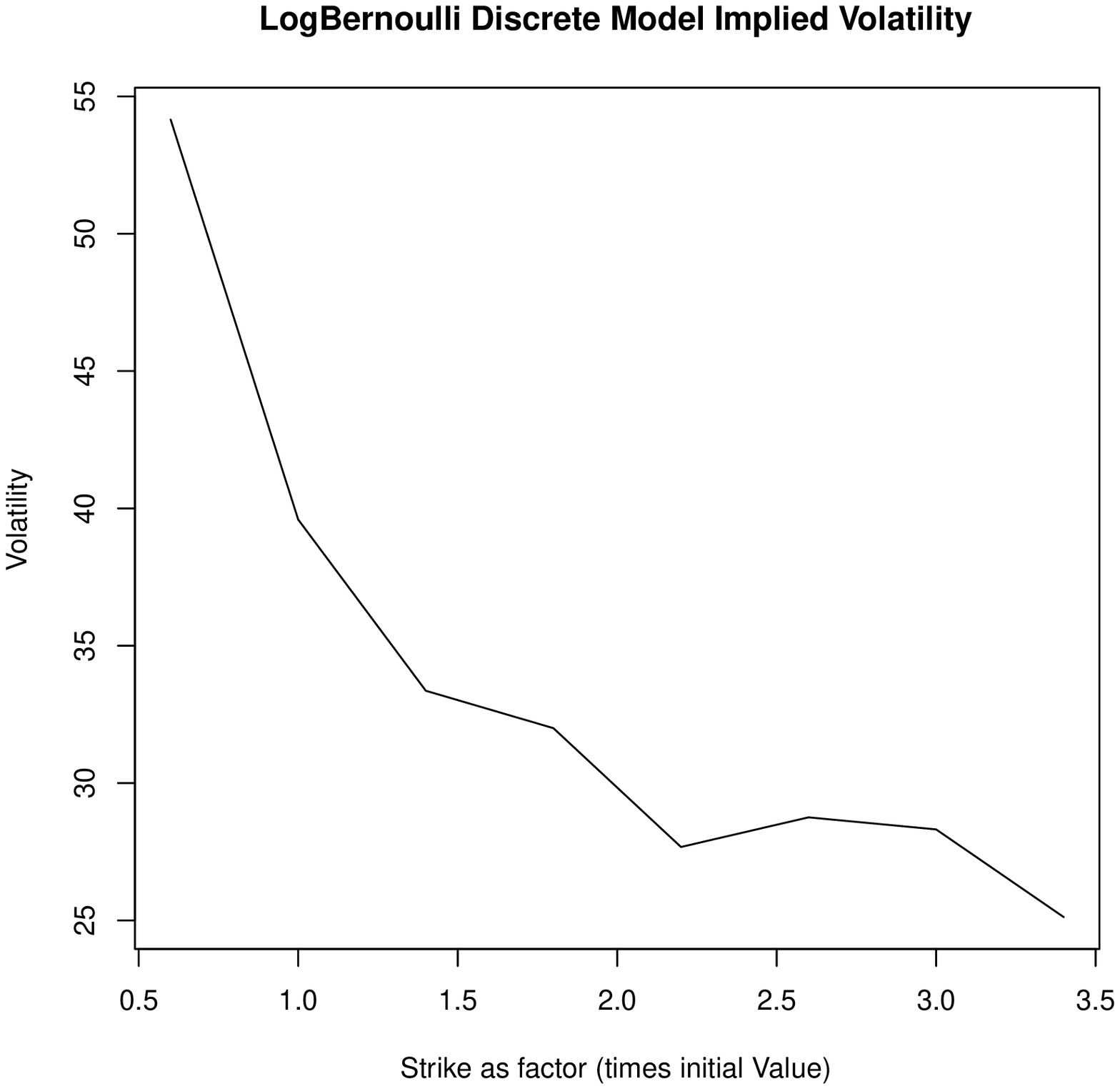}\includegraphics[width=4.8cm,height=4.5cm]{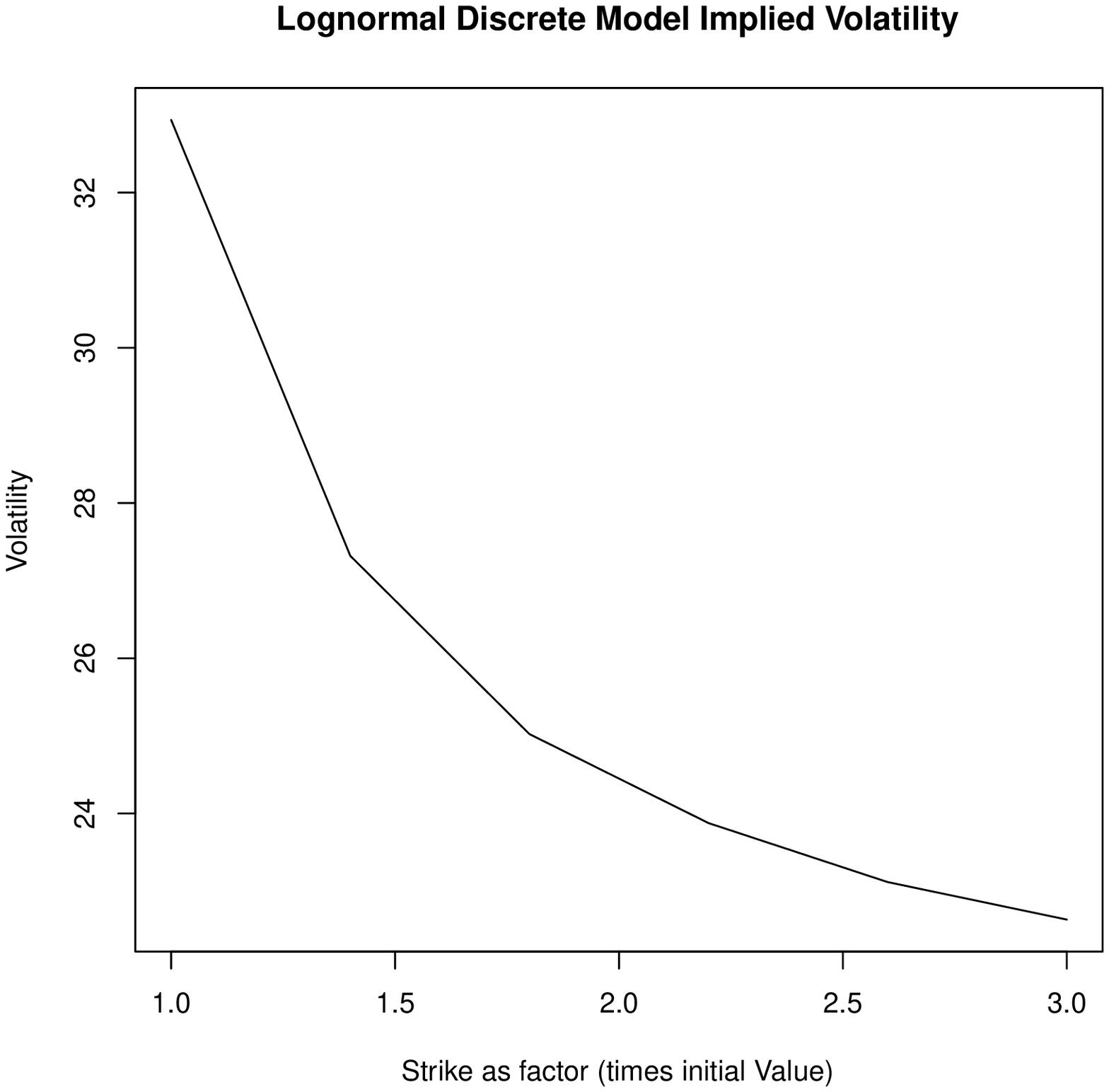}\includegraphics[width=4.8cm,height=4.5cm]{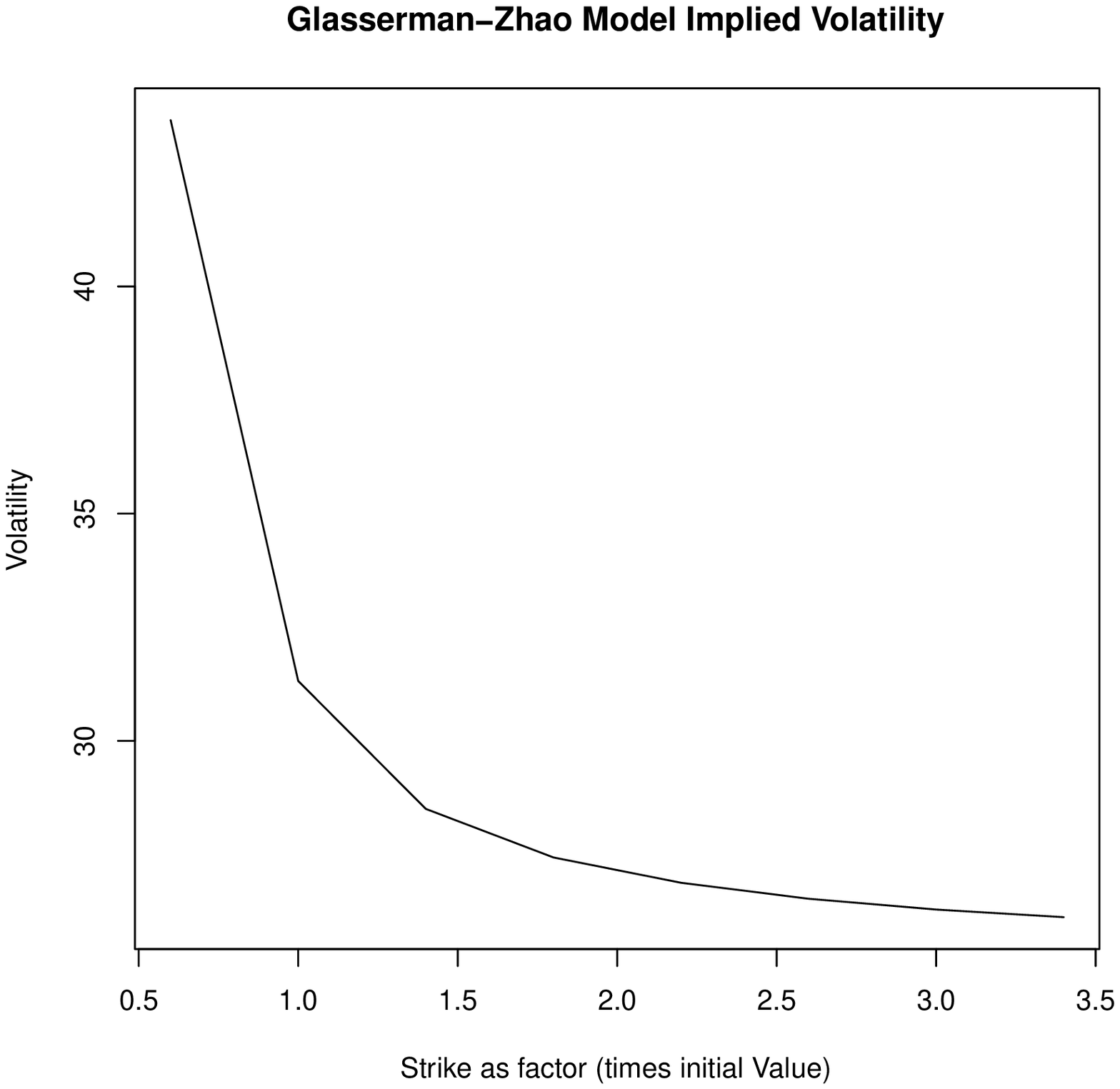}  \\
 If we look at the results from left to right, model $1$ (exact calculation with non frozen drift) yields considerably higher volatilities than the other $2$ models. Model $2$ (log normal discrete) yields higher implied volatilities than the Glasserman-Zhao discretization for caplets with small strikes but smaller volatilities as the strike grows. Model $3$ is one of the the Glasserman-Zhao arbitrage free discretizations from $\cite{GZ2000}$.\\
 The details on the implemented models and the convergence result for a refined time grid can be found in the following sections.
\section{Model Approach And Properties}\label{properties}
\subsection{Aims}
Based on what goals did we derive the discrete time LIBOR market analogue?
\begin{itemize}
\item
We expect the model to fit a given finite tenor structure with initial data for each point in the tenor structure.
\item
We expect the model to give an explicit expression for the evolution of the respective rates in the future.
\item
We expect the model to prevent arbitrage opportunities between different LIBOR rates.
\item
Finally, we want the model to approximate continuous time LIBOR market models if we let the discrete time grid become infinitely fine.
\end{itemize}
All of this can be accomplished using the following model:\\
\subsection{Discrete LIBOR Market Model}
We assume as given a stochastic basis $(\Omega,\mathcal F,\mathbf F,\mathbb P)$ with $\mathbf F=(\mathcal F_t)_{t\geq 0}$ satisfying the usual conditions.\\
Given a time horizon $T^*$, we represent running time as a finite grid of $m=(n+1)p$ positive real numbers ($n,p\in\mathbb N$)
\[
t_i:=\frac im T^*,\qquad i=0,1,\ldots ,m
\]
and we define an $(n+1)$-element tenor structure through
\[
T_j:=\frac j{n+1} T^*.\qquad j=1,\ldots, n+1.
\]
We further assume we have a family of positive real numbers $(\lambda_{ij})_{i=0,\ldots, m, j=1,\ldots, n}$ and positive real values $\{ L(0,T_j)\}_{j=0}^n$. Furthermore let $\delta :=\frac 1{n+1}$.
\begin{definition}[Discrete Forward LIBOR Market Model]\label{discLIB}
A family of stochastic processes $\{ L(t_i,T_j) \}_{i =0,\ldots ,pj ,j=1,\ldots, n}$ and a family of equivalent probability measures $\{\mathbb P_j \}_{j=2,\ldots,n+1}$ constitutes a LIBOR market model if
\begin{enumerate}
\item
$L(.,T_j)$ is a $\mathbb P_{j+1}$ martingale for $j=1,\ldots,n$.
\item
For $j=2,\ldots,n$
\begin{equation}
\frac{d\mathbb P_j}{d\mathbb P_{j+1}}(t_i)=\frac{1+ \delta L(t_i,T_j)}{1+\delta L(0,T_{j})}.
\end{equation}
\end{enumerate}
\end{definition}
Indeed such a model can be constructed by
\begin{theorem}[Existence Of A Discrete LIBOR Market Model]\label{exists}
Given an adapted time discrete process build from random variables $(X_{t_i})_{i=1,\ldots, m}$ such that 
\begin{equation}\label{integra}
\mathbb E_{\mathbb P_{n+1}}\Big{(} \exp (\sum_{i=1}^m\sum_{j=1}^n \lambda_{ij}X_{t_i})\Big{)}<\infty .
\end{equation}
There exists a discrete LIBOR market model $(L(.,T_j))_{j=1}^n$ 
\begin{equation}\label{LIBORdef}
L(t_i,T_{j})=L(0,T_j)  \exp{ ( \sum_{u=1}^i \lambda_{uj} (X_{t_u}+b_{t_u}^j) )},\qquad i=1,\ldots ,jp
\end{equation}
with
\begin{equation}\label{drift}
b_{t_i}^j= -\frac1{\lambda_{ij}}\log \mathbb E_{\mathbb P_{j+1}} ( \exp (\lambda_{ij}X_{t_i}) |\mathcal F_{t_{i-1}}).
\end{equation}
\end{theorem}
The proof is straight forward. It can be found explicitly in the appendix.\\
 This approach allows us to define implied Bond price and Forward dynamics:\\
 Let $\eta (t_i):=u$ for $ \frac{(u-1)}{n+1}T^*  < t_i \leq \frac u{n+1}T^*$.
\begin{definition}[Bondprices And Forward Price Process]\label{Bonfor}
We define the family of processes called bond price processes $(B(t_i,T_j))_{i=0,\ldots,pj, j=1,\ldots, n}$ from the LIBOR rate processes $\Big{(}L(.,T_j)\Big{)}$ through
\begin{equation}
B(T_i,T_j)=\prod_{u=i}^{j-1} \frac1{1+\delta L(T_i,T_u)},\qquad i\leq j-1
\end{equation}
and then
\begin{equation}\label{bond}
B(t_i,T_j)= B (0,T_{\eta (t_i)-1})\prod_{u=\eta (t_i)-1}^{j-1} \frac1{1+\delta L(t_i,T_u)}.
\quad \forall t_i \neq T_r, r=1,\ldots , n
\end{equation}
We call
\begin{equation}\label{forward}
F_B(t_i,T_j,T_{j+1} ):=1+\delta L(t_i,T_j), \quad T_j \in \{ T_1,\ldots, T^* \}, \delta = \frac 1{n+1}
\end{equation}
 the forward price process. $F_B(t_i,T,T+\delta)$ is a $\mathbb P_{T+\delta}$ martingale.
\end{definition}
With that, we get a full LIBOR Market Model.
\section{Convergence Result To Continuous Time LIBOR}\label{convergence}
While we view the discrete LIBOR market model analogue to be of intrinsic interest, the central idea when developing the model was to allow for an approximation of continuous time LIBOR market models without having to freeze the drift.\\
In this section we therefore derive a convergence result of a sequence of discrete time models to a continuous time model. First, we derive a difference equation representation to make the applicability of convergence results more immediate in the following.
\subsection{Difference Equation Presentation of the Dynamics}
We apply the following definitions to the process $\mathbf X_{t_i} + \mathbf B_{t_i}=  \sum_{u=1}^i (X_{t_u} + b_{t_u}) $.
\begin{definition}[Jump Measure]\label{Jump}
We call
\begin{equation} 
\mu (A,\{ t_i\}_{i=1}^j):= \#  \{ t_i\in \{ t_i \}_{i=0}^j | X_{t_i}+b_{t_i} \in A  \}\quad j>0
\end{equation}
 the jump measure of the process $\mathbf X_{t_i}+ \mathbf B_{t_i}$.
\end{definition}
With this definition it holds for $f$ continuous and $\mathbb E[f(X_{t_i})]<\infty, \forall i$
\begin{equation}
\int_{0}^{t_i}\int_{\mathbb R} f(x)\mu (dx,dt)= \sum_{j=1}^i f( X_{t_j}+b_{t_j}).
\end{equation}
Furthermore we need
\begin{definition}[Compensator]
We define the compensator of a random jump measure $\mu$ with respect to a measure $\mathbb P$ to be 
\begin{equation}
\nu (A,\{t_i \}_{i=1}^j)= \mathbb E_{\mathbb P} (\mu (A,\{ t_i \}_{i=1}^j)|\mathcal F_{t_{j-1}})
\end{equation}
\end{definition}

One can then show
\begin{theorem}[Dynamics And Measure Change]\label{Diffe}
Given a discrete LIBOR market model as constructed in theorem $(\ref{exists})$. Then the dynamics of the LIBOR rates $\Big{(}L(t_i,T_j)\Big{)}_{i=0,\ldots,jp , j=1,\ldots,n}$ under their respective forwardmeasures $(\mathbb P_j)_{j=2,\ldots,n+1}$ are
\begin{equation}
\Delta L(t_i,T_j)=L(t_{i-1},T_j)\int_{\mathbb R}(\exp (\lambda_{ij}x)-1)(\mu -\nu^{j+1} )(dx,t_i) \quad \textrm{under $\mathbb P_{j+1}$}
\end{equation}
where
\begin{equation}
\mu (A,\{t_i \}_{i=1}^j):= \#  \{ t_i\in \{t_i \}_{i=1}^j | X_{t_i}+b_{t_i} \in A  \}
\end{equation}
and 
\begin{equation}
\nu^{j+1} (A,\{t_i \}_{i=1}^u)= \mathbb E_{\mathbb P_{j+1}} (\mu (A,\{ t_i\}_{i=1}^u )|\mathcal F_{t_{u-1}}).
\end{equation} 
And then the measure change from $\mathbb P_{j+1}$ to $\mathbb P_j$ is reflected in a change of the compensator
\begin{equation}
\nu^{j} (dx,\{ t_i\} )=(\ell (t_{i-1},T_j)( \exp (\lambda_{ij}(x))-1)+1)\nu^{j+1}(dx,\{ t_i\} ).
\end{equation}
Here $\ell$ is completely analogous to continuous time:
\begin{equation}\label{ell}
\ell (t_{i-1},T_j)=\frac{\delta L(t_{i-1},T_j)}{1+\delta L(t_{i-1},T_j)}.
\end{equation}
\end{theorem}
Due to the difference equation representation, we can write the measure-change down as a change in the compensator, as in the continuous time case. 
\subsection{Approximation Theorem For Discrete LIBOR}\label{Theorem}
For the convergence result of the sequence of discrete time models, we embed the discrete models in a common space. This is best describes the way we think about the discrete time models if they are used for approximation of continuous time models.
\begin{theorem}[Approximation Theorem]\label{approximationtheorem}
Given a sequence $(X^{(k)})_{k\in\mathbb N}$ of discrete time PIIs on stochastic bases $\mathcal B^{(k)}=\Big{(}\Omega^{(k)},, \mathcal F^{(k)},(\mathcal F^{(k)}_i )_{i\in\mathbb N},\mathbb P_{n+1}^{(k)} \Big{)}$.\\
For each $k$ we define an adopted sequence of random variables 
\begin{equation}
U_i^k:=  X_i^{(k)}- X_{i-1}^{(k)},\quad i\in\mathbb N
\end{equation}
and a change of time $(\sigma_t^k)$ to consider then as processes with paths in $\mathbb D ([0,T^*]):=\{ f: [0,T^*]\rightarrow \mathbb R | f \textrm{is cadlag} \}$ with fixed jump times the processes 
\[
 X_t^{(k)}:=\sum_{1\leq i\leq \sigma_t^k} U_i^k,\qquad t\in \mathbb R_{+}.
\]
 As $k$ goes to infinity we let the mesh of the set of jump times go to $0$.\\
Now we define a PII $\Big{(} X \Big{)}$ on $\Big{(}\Omega,\mathcal F,(\mathcal F_t )_t, \mathbb P_{n+1} \Big{)}$
\begin{equation}
 X_t + B_t = \int_0^t b_s ds+\int_0^t \int_0^t c_s^{\frac12} dW_s +\int_0^t\int_{\mathbb R} x (\mu-\nu )(dx,ds),
\end{equation}
where $(\mathcal F_t)_{t\in [0,T^*]}$ is the filtration generated by $( X_t)_{t\in [0,T^*]}$ and
\begin{equation}
B_t = \int_0^t b_s ds.
\end{equation}
Also $(X_t)$ has to fulfill the conditions on a driving process for a L\'evy LIBOR model as in $\cite{EO2005}$. This imposes the following condition on the drift characteristic of the driving L\'evy process
\begin{equation}
\int_0^t \lambda (s,T^*-\delta ) b_s ds = \int_0^t c_s \frac 1 2 \lambda (s,T^*-\delta )^2 ds + \int_0^t \int_{\mathbb R} (e^{\lambda (s,T^*-\delta )x}-1-\lambda (s,T^*-\delta )x )v( ds, dx) .
\end{equation}
We assume that  
\begin{equation}
 X^{(k)}+ B^{(k)} \rightarrow  X + B ,\quad \textrm{weakly as processes}
\end{equation}
for
\begin{equation}
B^{(k)}_{t_i} = \sum_{j=1}^i b^{(k)}_{t_j}
\end{equation}
(where the $b^{(k)}$ are defined in $\ref{drift}$), 
and it holds that
\begin{equation}\label{condition}
\sup_{k\in\mathbb N \cup { \infty } }\mathbb E (\int_0^{T^*}\int_{|x|>1} \exp (ux)F^{(k)}_s(dx)ds)<\infty, u\leq (1+\epsilon )M, M\geq \sum_{j=1}^n |\lambda (.,T_j)|,  
\end{equation}
\[
 \nu^{(k)} (dx,ds) =F^{(k)}_s(dx)ds., \nu^{\infty}(dx,ds)=F^{\infty}_s(dx)ds=\nu (dx,ds)=F_s(dx)ds.
\]
as well as
\begin{equation}
\sup_{k\in\mathbb N\cup{ +  \infty }} \mathbb E(\int_0^{T^*}\int_{\mathbb R} (x^2 \land 1)F^{(k)}_s(dx)ds)<\infty .
\end{equation}
We define discrete LIBOR market models $(L(t_i,T_j)^{(k)})_{i=0,1,\ldots ,jp^k, j=1,\ldots ,n,k\in\mathbb N}$, driven by the $(  X^{(k)})_{k\in\mathbb N}$ as in the section on discrete forward LIBOR modeling. In other words the rate $L(.,T_j)^{(k)}$ is driven by the exponential of $( X^{(k)})$ minus a drift, making the exponential of $(X^{(k)})$ into a martingale under $\mathbb P_{j+1}^{(k)}$.\\
We assume that for $i^{k}:= \sup_{u =0,1,\ldots,jp^k } \frac {u}{np^k}\leq t_i $ it holds that $\lambda_{i^{k}j}^{k}\rightarrow \lambda(t_i,T_j)$ pointwise, and $\lambda(.,T_j):\mathbb R_+\rightarrow \mathbb R_+$ are linearly bounded, locally Lipschitz and positive functions. Furthermore we assume we are given starting values $\Big{(}L(0,T_j)\Big{)}_{j=1,\ldots ,n}$ which are strictly increasing in $j$ then 
\begin{equation}
L(.,T_j)^{(k)}\textrm{under $\mathbb P^{(k)}_{j+1}$}\rightarrow L(.,T_j)\quad\textrm{under $\mathbb P_{j+1}$}\qquad \forall 1\leq j\leq n
\end{equation}  
weakly as a process.
 \end{theorem}

\begin{lemma}[Contiguity of forward measures]\label{contig}
Under the conditions of the preceding theorem $(\ref{approximationtheorem})$, any two sequences of forward measures are contiguous. 
\end{lemma}
\section{Implementation Results}\label{implementation}
In order to support the practical use of the discrete time LIBOR market model, we calculated caplet prices and volatility smiles for a reference rate. This shows that not only is pricing of options feasible based on the discrete time model, but the results allow to find other arbitrage free prices for discrete approximation of the lognormal LIBOR market model.\\
The discrete model $(\ref{discLIB})$ was implemented using
\begin{enumerate}
\item
A purely discrete driving process with Bernoulli distributed random variables $(X_{t_i})$ with parameter $p$. For this model, it was well possible to calculate the exact LIBOR model. 
\item
A model driven by standard normal $(X_{t_i})$ . This model was used for Monte Carlo simulation, as the exact local drift would be to costly to calculate.
\end{enumerate}
As a reference, an implementation of one of the models described in $\cite{GZ2000}$, using the same volatilities $(\lambda_{ij})$, was carried out.\\
\subsection{General Setup}
Our aim was to calculate the price of forward caplets on $L(T_5,T_5)$ in a model with $10$ LIBOR rates. Hence the tenor structure was 
\begin{equation}
T_i = i. \quad i = 1,\ldots, 11
\end{equation}
with the $T_i$ denoting years and the step size being $\delta = 1$ giving us a time grid $t=0,1,2,\ldots, 11$. The  time $T_{11}$ was not used to model a LIBOR rate but only to define the terminal measure $\mathbb P_{T_{11}}$ and hence the forward measure for $(L(t_i,T_{10}))_{t_i\in 0,\ldots, 10}$.\\
Bond prices $B(0,T_2),\ldots, B(0,T_{11})$ were calculated from $(\ref{bond})$.\\
The caplet price to be calculated was
\begin{equation}
\Pi = B(0,T_6)\mathbb E_{\mathbb P_{6}} [L(T_5,T_5)- KL(0,T_5) ]
\end{equation}
for $K$ values 
\\
\vspace{0.2cm}
\begin{tabular}{|c|c|c|c|c|c|c|c|}
\hline
   $K_1$  & $K_2$ & $K_3$  & $K_4$ & $K_5$ & $K_6$ & $K_7$  &  $K_8$   \\
  $0.6$ & $1$ & $1.4$  & $1.8$ & $2.2$ & $2.6$ & $3$ & $3.4$ \\
 \hline
\end{tabular}\\
\vspace{0.2 cm}
The initial data was\\
\vspace{0.2cm}
\begin{tabular}{|c|c|c|c|c|c|c|c|c|c|}
\hline
   $L(0,T_1)$  & $L(0,T_2)$ & $L(0,T_3)$  & $L(0,T_4)$ & $L(0,T_5)$ & $L(0,T_6)$ & $L(0,T_7)$  &  $L(0,T_8)$  & $L(0,T_9)$  & $L(0,T_{10})$ \\
  $0.0207$ & $0.23$ & $0.0262$  & $0.28$ & $0.0292$ & $0.0318$ & $0.0342$ & $0.0362$ & $0.0379$ & $0.04$\\
 \hline
\end{tabular}\\
\vspace{0.2 cm}
For comparison purposes, all chosen models have driving process variables $X_{t_i}$ with $\mathbb E[X_{t_i}]=0$ and $\mathbb V[X_{t_i}]=1$. All models converge to the lognormal LIBOR market model for refined time grids.\\
 As a measure of how the model results related to one another, we compared the caplet prices, in terms of implied volatilities, obtained for the below a priori volatility structure for all $3$ models.\\
The a priori volatility parameters were  \\
\vspace{0.2 cm}
\begin{tabular}{|c|c|c|c|c|c|c|c|c|c|c|}
\hline
  Volatilities $(\lambda_{ij})$  &  $\lambda_{i1}$  & $\lambda_{i2}$ & $\lambda_{i3}$  & $\lambda_{i4}$ & $\lambda_{i5}$ & $\lambda_{i6}$ & $\lambda_{i7}$  &  $\lambda_{i8}$  & $\lambda_{i9}$  & $\lambda_{i10}$ \\
  \hline
 & $0.34$ & $0.32$ & $0.3$  & $0.28$ & $0.26$ & $0.24$ & $0.22$ & $0.2$ & $0.18$ & $0.16$\\
 \hline
\end{tabular}\\
\vspace{0.2 cm}
The volatility parameters were independent of the time index $i$ and depended on the rate index $j$.\\
\subsection{Implementation Of The Log-Bernoulli Model}
Based on the central equation from section $(\ref{properties})$, Theorem $(\ref{exists})$, equation $(\ref{LIBORdef})$
\begin{equation}
L(t_i,T_{j})=L(0,T_j)  \exp{ ( \sum_{u=1}^i \lambda_{uj} (X_{t_u}+b_{t_u}^j) )}
\end{equation}
using Bernoulli Variables $(X_{t_i})$, $\mathbb P [X_{t_i}=1]=0.5$ and $\mathbb P [X_{t_i}=-1]=0.5$, it was possible to calculate all paths of the resulting discrete model. The pricing of a caplet then reduces to summing up over all possible payoff values times the respective path probability.\\
 While the number of paths is $2^5$, in this example, and we have to calculate from the terminal rate to the target rate $L(,5)$, the speed of the computation is fast for the sample data.\\
 The main computational issue is to properly calculate
\begin{equation}
b_{t_i}^j= -\frac1{\lambda_{ij}}\log \mathbb E_{\mathbb P_{j+1}} ( \exp (\lambda_{ij}X_{t_i}) |\mathcal F_{t_{i-1}})=
\end{equation}
\[
-\frac1{\lambda_{ij}}\log \mathbb E_{\mathbb P_{n+1}} ( \exp (\lambda_{ij}X_{t_i}) \prod_{k=j+1}^{n}\frac{F_B(t_i,T_k,T_{k+1})}{F_B(0,T_k,T_{k+1})}|\mathcal F_{t_{i-1}}).
\]
Differently put, we need the compensator $v^{j}$ under the terminal measure $\mathbb P_{n+1}$
\begin{equation}
v^{j} (dx,\{ t_i\} )=(\prod_{k=j}^n (\ell (t_{i-1},T_k)( \exp (\lambda_{ik}(x))-1)+1)))\nu^{n+1}(dx,\{ t_i\} ).
\end{equation}
We therefore need a different compensator value depending on time and the $\ell (t_{i-1},T_k)$ of the earlier rates. For a model with $10$ rates this proved to be no issue, once the compensator product was implemented. Our implementation is based on the fact that
\begin{equation}
(\prod_{k=j}^n (\ell (t_{i-1},T_k)( \exp (\lambda_{ik}(x))-1)+1)))\nu^{n+1}(dx,\{ t_i\} ) = 
\end{equation}
\[
\sum_{\sigma}(\prod_{k\in\sigma} (\ell (t_{i-1},T_k)( \exp (\lambda_{ik}(x))))\prod_{k\in \sigma^C} (1-\ell (t_{i-1},T_k))).
\]
where $\sigma$ denotes all subsets of the set $\{ j, j+1,\ldots, n \}$ and $\sigma^C$ is the complement of $\sigma$ in $\{ j, j+1,\ldots, n \}$.\\
The resulting implied volatilities are
\\
\vspace{0.2cm}
\begin{tabular}{|c|c|c|c|c|c|c|c|c|}
\hline
  K & $0.6$ & $1$ & $1.4$  & $1.8$ & $2.2$ & $2.6$ & $3$ & $3.4$ \\
\hline
  implied Volatility  & $0.542$ & $0.396$  & $0.334$ & $0.32$ & $0.277$ & $0.2875$ & $0.2832$ & $0.25$   \\
 \hline
\end{tabular}\\
\vspace{0.2 cm}
In the figure below, the implied volatilities are rescaled by a factor of $100$.
\begin{center}
Volatility Smile for the exact discrete time log-Bernoulli model \\
\includegraphics[width =7 cm,height= 7 cm]{LogBernoulli.eps}
\end{center}
\subsection{Implementation Of The Lognormal Discrete Analogue}
We again used section $(\ref{properties})$,theorem $(\ref{exists})$ ,equation $(\ref{LIBORdef})$
\begin{equation}
L(t_i,T_{j})=L(0,T_j)  \exp{ ( \sum_{u=1}^i \lambda_{uj} (X_{t_u}+b_{t_u}^j) )}
\end{equation}
but here with $X_{t_u}$ being standard normally distributed. The compensator under the terminal measure is once more
\begin{equation}
v^{j} (dx,\{ t_i\} )=(\prod_{k=j}^n (\ell (t_{i-1},T_k)( \exp (\lambda_{ik}(x))-1)+1)))\nu^{n+1}(dx,\{ t_i\} )
\end{equation}
but here $\ell (t_{i-1},T_k)$ can have infinitely many possible values and the path by path calculation of the above model becomes impossible. Therefore, pricing was carried out, using $500 000$ paths from Monte Carlo Simulations. This model is somewhat slower than the Glasserman-Zhao discretization in the present implementation, due to the compensator product. It however yields different results and therefore yields other arbitrage free prices.
The resulting implied volatilities are
\\
\vspace{0.2cm}
\begin{tabular}{|c|c|c|c|c|c|c|c|c|}
\hline
  K & $0.6$ & $1$ & $1.4$  & $1.8$ & $2.2$ & $2.6$ & $3$ & $3.4$ \\
\hline
  implied Volatility  & $0.518$ & $0.333$  & $0.276$ & $0.253$ & $0.241$ & $0.234$ & $0.23$ & $0.226$ \\
 \hline
\end{tabular}\\
\vspace{0.2 cm}
In the figure below, the implied volatilities are rescaled by a factor of $100$.
\begin{center}
Volatility Smile for the lognormal disrete model - Monte Carlo Simulations\\
\includegraphics[width =7 cm,height= 7 cm]{Lognormalsmile.eps}
\end{center}
\subsection{Implementation Of One Glasserman-Zhao Arbitrage Free Discretization}
 For the arbitrage free discretization of $\cite{GZ2000}, p.43$ we used the equations $(21)$, $(24)$ and $(25)$ with $\delta=1$ and $h=1$
\begin{equation}
L(t_i,T_j) = \frac{W_j(i)}{1+ W_{j+1}(i)+\ldots + W_{n}(i)} 
\end{equation}
with the same $i,j$ as in the above models except for $(L(t_i,T_10))_{t_i \in 0, \ldots ,10}$ which is modeled arbitrage free already. Here
\begin{equation}
W_j(i)=W_j(i-1)\exp \Bigg{(} -\frac12 \sigma_j^2 (i) + \sigma_j (i) Y_i  \Bigg{)}
\end{equation}
with
\begin{equation}
\sigma_j (i) =\lambda_{ij} +\sum_{k=j+1}^n \frac{ W_k \lambda_{ik}}{1 + W_k+\ldots +W_n}
\end{equation}
and $Y_i\sim\mathcal N (0,1)$ and the $(Y_i)_{i=1,\ldots, 10}$ are independent.\\
The pricing formula we used was from $\cite{GZ2000}$, on page $45$ adapted for $N=10$ and $n=5$.\\
 Here $\sigma_j(i)$ again includes local terms, dependent on the value for a higher maturity $W_{j+1}$. Therefore, for a reference implementation, we also used Monte Carlo Simulation for pricing ($500 000$ paths). There is no compensator product here and compared to the log normal discrete model, this method is faster. The log-Bernoulli model however turns out to be even faster (even using exact calculations) than this method. Test for further data will follow in future work.
 \vspace{0.2cm}
\begin{tabular}{|c|c|c|c|c|c|c|c|c|}
\hline
  K & $0.6$ & $1$ & $1.4$  & $1.8$ & $2.2$ & $2.6$ & $3$ & $3.4$ \\
\hline
  implied Volatility  & $0.437$ & $0.313$  & $0.285$ & $0.2745$ & $0.269$ & $0.265$ & $0.263$ & $0.261$ \\
 \hline
\end{tabular}\\
\vspace{0.2 cm}
In the figure below, the implied volatilities are rescaled by a factor of $100$.
\begin{center}
Volatility Smile for one arbitrage free discretization of Glasserman and Zhao \\
\includegraphics[width =7 cm,height= 7 cm]{GLZp.eps}
\end{center}
\section{Summary And Outlook}
We have proven the existence of the arbitrage free discrete LIBOR analogue $(\ref{discLIB})$ and its convergence to a continuous time model under reasonable assumptions (convergence of the driving processes, uniform integrability).\\
Further work will be along the lines of the work by Musiela and Rutkowski $\cite{MR97}$. We aim to establish a complete discrete LIBOR framework, namely forward models and the implied savings account. Also we want to adapt the convergence result $(\ref{approximationtheorem})$ for joint convergence under terminal measure.\\
Establishing some relation to the general HJM framework would be exciting.\\
Interest in LIBOR Market Models is high as ever as the work of Keller-Ressel, Papapantoleon and Teichmann $\cite{label6084}$ and recent work by Eberlein et al show.
\newpage
\appendix
\section{Proofs}\label{proofs}
\textbf{Proof of theorem $\ref{exists}$, section $(\ref{properties})$}:\\
First we check what $L(.,T_j)$ has to fulfill to be a $\mathbb P_{j+1}$ martingale.\\
We proceed inductively: Assume $\Big{(}L(t_{i},T_j)\Big{)}_{i=0,1,\ldots ,u-1}$ has been constructed to be a martingale. We look at $L(t_u,T_j)$
\begin{equation}
\mathbb E_{\mathbb P_{j+1}} ( L(t_u,T_j) |\mathcal F_{t_{u-1}})=L(t_{u-1},T_j)\mathbb E_{\mathbb P_{j+1}} ( \exp (\lambda_{uj}(X_{t_u}+b_{t_u}^j)) |\mathcal F_{t_{u-1}}).
\end{equation}
We therefore have as martingale condition
\begin{equation}
\mathbb E_{\mathbb P_{j+1}} ( \exp (\lambda_{uj}(X_{t_u}+b_{t_u}^j)) |\mathcal F_{t_{u-1}})=1\Leftrightarrow \mathbb E_{\mathbb P_{j+1}} ( \exp (\lambda_{uj} X_{t_u} ) |\mathcal F_{t_{u-1}})=\exp (-\lambda_{uj} b_{t_u}^j).
\end{equation}
Here we choose $b_{t_u}^j$ to be $\mathcal F_{t_{u-1}}$-measureable.\\
 This yields
\begin{equation}
b^j_{t_u}=-\frac1{\lambda_{uj}}\log \mathbb E_{\mathbb P_{j+1}} ( \exp (\lambda_{uj}X_{t_u}) |\mathcal F_{t_{u-1}}).
\end{equation}
So we have that $L(t_u,T_j)$ fulfills the martingale condition with this choice of drift.\\
Now assume we are given $(L(t_i,T_{j}))_{i=0,\ldots, nj}$ under $\mathbb P_{j+1}$ as a $\mathbb P_{j+1}$ martingale. We define $\mathbb P_{j}$ through
\begin{equation}
\frac{d\mathbb P_{j}}{d\mathbb P_{j+1}}(t_i):= \frac{1+\delta L(t_{i},T_{j})}{1+\delta L(0,T_{j})}.
\end{equation}
All appearing expectations exist due to the integrability condition $(\ref{integra})$. $\qquad \Box$\\
\vspace{0.4 cm}
\textbf{Proof of theorem $(\ref{Diffe})$, section $(\ref{convergence})$}:\\
We have that
\begin{equation}
\Delta L(t_i,T_j)=L(0,T_j)\exp \Big{(}\sum_{u=1}^{i-1} \lambda_{uj}(X_{t_u}+b_{t_u}^j) \Big{)}(\exp (\lambda_{ij}(X_{t_i}+b_{t_i}^j)) -1) =
\end{equation}
\[
L(t_{i-1},T_j)(\exp (\lambda_{ij}(X_{t_i}+b_{t_i}^j)) -1).
\]
And with the definition of $\mu$ and $\nu^{j+1}$ it holds that
\begin{equation}
\int_{\mathbb{R}}  (\exp{(\lambda_{ij}(x))}-1) \mu (dx, \{ t_i \})=
\exp{(\lambda_{ij}(X_{t_i}+b_{t_i}^j))}-1
\end{equation}
and
\begin{equation}
-\int_{\mathbb{R}}  (\exp{(\lambda_{ij}(x))}-1) \nu^{j+1}(dx,\{ t_i\})=
-\exp{(0)}+1=0.
\end{equation}
So  
\begin{equation}\label{Difference}
\Delta L(t_i,T_j)=L(t_{i-1},T_j)\int_{\mathbb{R}} (\exp{(\lambda_{ij}(x))}-1) (\mu-\nu^{j+1})(dx,\{ t_i\}).
\end{equation}
From the definition of $\nu^{j+1}$ it follows that
\begin{equation}
\int_{\mathbb{R}}  (\exp{(\lambda_{ij}x)} \nu^{j+1}(dx,\{ t_i\})= 1.
\end{equation}
We then proceed to look at the forward price process $F(t_i,T_{j},T_{j+1})=1+\delta  L(t_i,T_{j})$. The dynamics under $\mathbb P_{j+1}$ are derived as 
\begin{equation}
\Delta F(t_i,T_{j},T_{j+1})=\delta \Delta L(t_i,T_{j}).
\end{equation}
We have described $\Delta L(t_i,T_{j})$ above. Therefore we get an analogy to the SDE of the forward price in continuous time
\begin{equation}
\Delta F(t_i,T_{j},T_{j+1})=F(t_{i-1},T_{j},T_{j+1})\int_{\mathbb{R}} 
\ell (t_{i-1},T_{j})(\exp{(\lambda_{ij}(x)}-1) (\mu-\nu^{j+1})(dx,\{ t_i\}).
\end{equation}
From the proof of theorem $(\ref{exists})$ we know that
\begin{equation}
\frac{d\mathbb P_{j}}{d\mathbb P_{j+1}}(t_i)=\frac{F_B(t_i,T_j,T_{j+1})}{F_B(0,T_j,T_{j+1})}.
\end{equation}
 With the discrete Girsanov Theorem and the stochastic difference equation representation of $F_B(t_i,T_j,T_{j+1})$, we get
\begin{equation}
\frac{d\mathbb P_{j}}{d\mathbb P_{j+1}}(t_i)=Z^j_i=\mathcal E (M^j_i)
\end{equation}
with
\begin{equation}
M^j_i=\sum_{u=1}^i\int_{\mathbb R}\ell (t_{u-1},T_j)(e^{\lambda_{uj}(x)}-1)(\mu-\nu^{j+1})(dx,\{ t_u\} ).
\end{equation}
We apply discrete Girsanov theorems to get $Y_i^j$
\begin{equation}
Y_i^j(X_i)=\int_{\mathbb R}\Big{(}\ell (t_{i-1},T_j)(e^{\lambda_{ij}(x)}-1)+1\Big{)}(\mu-\nu^{j+1})(dx,\{ t_i\}).
\end{equation}
Therefore 
\begin{equation}
y_i(x_1,\ldots, x_{i-1};x_i)=\ell (t_{i-1},T_j)(e^{\lambda_{ij}(x_i)}-1)+1.
\end{equation}
We have that 
\begin{equation}
df^j_i(x_1,\ldots,x_{i-1};x_i)=\nu^j (dx,\{ t_i\})\quad \land \quad df^{j+1}_i(x_1,\ldots ,x_{i-1};x_i)=\nu^{j+1} (dx,\{ t_i\})
\end{equation}
and
\begin{equation}
y_i(x_1,\ldots,x_{i-1};x_i)=\frac{df^j_i}{df^{j+1}_i}(x_1,\ldots ,x_{i-1};x_i).
\end{equation}
 This gives us 
\begin{equation}
\nu^{j}(dx,\{ t_i\})=(\ell (t_{i-1},T_{j})(\exp{(\lambda_{ij}(x))}-1)+1)\nu^{j+1}(dx,\{ t_i\}).
\end{equation}
And this concludes the proof. $\quad\Box$\\
\vspace{0.4 cm}
\textbf{Proof of theorem $(\ref{approximationtheorem})$, section $(\ref{convergence}) $}:\\
We denote by
\begin{equation}
t_i^k:= \inf  \{ t \in \mathbb R_+|\sigma^k_t = i \}
\end{equation}
the embedded grid points of the discrete time models. Taking limits from the righthand side, we get weak convergence everywhere, after we have proved the result for the grid points.\\
With condition
\begin{equation}
\sup_{k\in\mathbb N\cup{ +  \infty } }\mathbb E (\int_0^{T^*}\int_{|x|>1} \exp (ux)F^{(k)}_s(dx)ds)<\infty, u\leq (1+\epsilon )M, M\geq \sum_{i=1}^n |\lambda (.,T_i)|,  
\end{equation}
 fulfilled, there exists a continuous time LIBOR model $L(t,T_n)$ with driving process $X$ as in $\cite{EO2005}$ and discrete Models $(L(t_i,T_n)^{(k)})_{k \in \mathbb{N}}$ with increments $U_i^{n,(k)}$.\\
The functions $\lambda(.,T_j)$ are linearly bounded and locally Lipschitz and the $(\lambda_{ij}^{(k)})_{i=0,\ldots ,p^{(k)}j, j=1,\ldots ,n}$ are pointwise convergent, deterministic functions, which we need for $(\ref{Approx})$ in Jacod and Shirjaev (\cite{JS2002}, page $578$, theorem $6.9$).
We have the stochastic difference equation,
\begin{equation}
\Delta L^{(k)}(t_i,T_{n})= L^{(k)}(t_{i-1},T_{n})\int_{\mathbb{R}} (\exp{(\lambda_{in}(x+b_{t_i}^{n,(k)}))}-1) (\mu -\nu^{n+1,(k)})(dx,t_i).
\end{equation}
To see that this converges to the SDE of $L(t,T_n)$ we need conditions on the convergence of stochastic differential equations in a general semimartingale setting. Such conditions can be found for instance in Jacod,Shiryaev $\cite{JS2002}$, page $577$f . \\
First we need certain existence and uniqueness results:\\
In our concrete case and using the notation from the appendix, we have that:
\[
Y=Z+f(Y_-)\cdot O
\]
with
\[
Y^{(k)}_{t_i}=L^{(k)}(t_i,T_n)\qquad \forall k \in\mathbb{N}
\]
\[
Z^{(k)}=L^{(k)}(0,T_n)\qquad \forall k \in\mathbb{N}
\]
\[
f^{(k)}(x)= x \qquad \forall k \in\mathbb{N}
\]
\[
O^{(k)}_{t_i}=\int_{\mathbb R} (\exp{(\lambda_{in}(x))}-1) (\mu -\nu^{n+1,(k)})(dx,t_i)\qquad \forall k \in\mathbb{N}.
\]
and
\begin{itemize}
\item
Our $ O^{(k)}$ is weakly convergent by assumption to
\begin{equation}
O := \int_{\mathbb R} (\exp (\lambda (t,T_{n+1})x)-1)(\mu -\nu^{n+1})(dx,dt).
\end{equation}
 Therefore certainly tight. In absence of a drift component, we automatically have that the drifts variation norm $Var(B^{n,(k)})_t$ (see \cite{JS2002}, page $27$, §$3a.$) is tight. Therefore $(O^{(k)})_k$ is $P-UT$ (see \cite{JS2002}, page $380$, $6.15$ and page $382$, $6.21$ or appendix $B$, definition $(\ref{TPUT})$). For the definition of $P-UT$ see \cite{JS2002}, page $377$, $6.1$ or appendix $B$, definition $(\ref{PUT})$.
\item
We have that $f^{(k)}$ is Lipschitz and linearly bounded for all $(k)$.
\item
We have assumed $O^{(k)} \rightarrow O$ weakly and shown that it is also $P-UT$ and know that $Z^{(k)}=L(0,T_n)\rightarrow L(0,T_n)$. Therefore we have $(O^{(k)},Z^{(k)})\rightarrow (O,Z)$ weakly and by the theorem on weak convergence of solutions of SDE's we have $(O^{(k)},Z^{(k)},Y^{(k)})\rightarrow (O,Z,Y)$ weakly.
\end{itemize}
Those conditions fulfilled entail that the sequence of discrete processes converges weakly to a weak solution of the proper SDE in continuous time.\\
That SDE is
\begin{equation}
dL(t,T_n)=L(t_{-},T_n)(c_s^{\frac12} dW_s+\int_{\mathbb R}(e^{\lambda (t,T_n)x}-1)(\mu -\nu^{n+1,(k)})(dx,dt))
\end{equation}
By the same reasoning, we get that $\Delta F_B( t_i,T_{n},T_{n+1})$ also converges to the SDE of the forward in continuous time.\\
By contiguity of the forward measures (see section $(\ref{convergence})$,lemma $(\ref{contig})$ and the proof below) all further rates converge under their respective forward measures. $\qquad\Box$\\
\vspace{0.4 cm}
\textbf{Proof of lemma $(\ref{contig})$, section $(\ref{convergence})$}:\\
We denote by
\begin{equation}
t_i^k:= \inf \{ t\in \mathbb R_+ | \sigma^k_t = i \}
\end{equation}
the embedded grid points of the discrete models and prove the result for the grid points.\\
The measure change between two consecutive forward measures is given by 
\begin{equation}
\frac{d\mathbb P_{j}}{d\mathbb P_{j+1}}(t)=Z^j_{t}=\mathcal E (M^j_t)=\frac{F_B(t,T_j,T_{j+1})}{F_B(0,T_j,T_{j+1})}
\end{equation}
and with respect to the terminal measure
\begin{equation}
\frac{d\mathbb P_{j}}{d\mathbb P_{n+1}}(t)=\prod_{l=j}^{n} \frac{F_B(t,T_l,T_{l+1})}{F_B(0,T_l,T_{l+1})}.
\end{equation}
Let $Z^{j,(k)}$ denote the measure change from $\mathbb P^{(k)}_{j+1}$ to $\mathbb P^{(k)}_{n+1}$.\\ 
With this notation we use $\cite{JS2002}$ p.$288$, theorem $1.10, (ii')$. Our sequence is uniformly integrable due to $(\ref{condition})$. This yields contiguity according to $\cite{JS2002}$ p.$288$, lemma $1.10, (ii')$. $\quad\Box$

\section{Weak convergence of processes}\label{weak}
\begin{definition}[Predictably Uniformly Tight]\label{PUT}
For each integer $n$ let $\mathcal{B}^n=(\Omega^n, \mathcal{F}^n,\mathbf{F}^n,\mathbb{P}^n)$ be a stochastic basis. We denote by $\mathcal{H}$ the set of all predictable processes $H^n$ on $\mathcal{B}^n$ having the form
\[
H^n_t=Y_01_{(0)}+ \sum_{i=1}^k Y_i 1_{(s_i,s_{i+1}]}(t)
\]
with $k\in \mathbb{N}$, $0=s_0<s_1<...<s_k<s_{k+1}$ and $Y_i$ is $\mathcal{F}^n_{s_i}$-measurable with $|Y_i|\leq 1$.\\
We can define an elementary stochastic integral by
\[
H^n \cdot  X^n_t=\sum_{i=1}^k Y_i(X_{\inf{\{ t,s_{i+1}\}}}^n-X_{\inf{\{ t,s_{i}\}}}^n).
\]
Now a sequence $(X^n)$ of adapted( on their respective Bases) cadlag $d$-dimensional processes is P-UT if for every $t>0$ the family of random variables $(\sum_{1\leq i \leq d}H^{n,i}\cdot X^{n,i}_t : n\in\mathbb{N},H^{n,i}\in\mathcal{H}^n)$ is tight in $\mathbb{R}$ meaning
\[
\lim_{a\rightarrow \infty}\sup_{H^{n,i}\in\mathcal{H}} \mathbb{P} (|\sum_{i=1}^d H^{n,i}\cdot X_t^{n,i}|>a)=0.
\]
\end{definition}
Assume the following setting:\\
\begin{itemize}
\item
We are given an equation
\[
Y=Z+f(Y_-)\cdot O
\]
as above
\item
For each n we have a stochastic basis $\mathcal{B}^n=(\Omega^n, \mathcal{F}^n, \mathbf{F}^n, \mathbb{P}^n)$ and an equation
\begin{equation}
Y^n=Z^n+ f_n(Y^n_-) \cdot  O^n.
\end{equation}
\item
$O^{n}$ is a d-dimensional-semimartingale on that basis.
\item
$Z^{n}$ is a q-dimensional cadlag adapted process.
\item 
$f_n$ are functions $\mathbb{R}^q\rightarrow \mathbb{R}^q \times \mathbb{R}^d$ such that each equation above admits a unique solution.
\end{itemize}
Then there holds
\begin{theorem}[Weak Approximation]\label{Approx}
Assume the functions $f_n$ fulfill Lipschitz and linear boundedness with constants not dependent on $n$ and $f_n\rightarrow f$ at least pointwise. Assume further that the sequence $O^n$ is P-UT. Then if $ Y^n$ denotes the unique solution of the sequence of equations there holds\\
If $(O^n,Z^n)\rightarrow (O,Z)$ weakly, then $(O^n,Z^n,Y^n)\rightarrow (O,Z,Y)$ weakly.
\end{theorem}
Proof:\\
See Jacod and Shirjaev \cite{JS2002}, page $578$, theorem $6.9$. $\qquad \Box$\\
\\
In the case of the theorem on discrete approximation we need to show especially $P-UT$ of the driving process of the SDE. However we have weak convergence of the sequence of driving processes already and in that case there holds
\begin{theorem}[Tightness And P-UT]\label{TPUT}
Let $(O^n)_n$ be a sequence of d-dimensional semimartingales with characteristics and second modified characteristics $(B^n,C^n,\nu^n)$ and $\tilde{C}^n$.\\
If the sequence $(O^n)_n$ is tight, then the sequence is $P-UT$ iff the sequence $Var(B^{n,i})_t$ is tight.
\end{theorem}
Proof:\\
Using Jacod and Shirjaevs book \cite{JS2002} we combine page $380$, $6.15$ and page $382$, $6.21$ and the theorem follows.  $\qquad \Box$\\
\newpage

\bibliography{dissertatio}

\vfill\tt\jobname.tex~(\today)

\end{document}